\documentclass[12 pt]{amsart}

\usepackage{amscd,amsfonts,amssymb,amsmath,amsthm,enumerate,latexsym,floatflt,longtable}
\usepackage{xcolor}
\usepackage{placeins}
\usepackage{hyperref}
\usepackage[utf8]{inputenc}

\title[The Table of the Structure Constants]{The Table of the Structure Constants for the Complex Simple Lie Algebra of Type $G_2$
and Chevalley Commutator Formulas in the Chevalley Group of Type $G_2$ over a Field}

\author[Sergey G.~Kolesnikov]{Sergey G.~Kolesnikov}

\date{\today}

\begin{document}
\sloppy

\maketitle
\begin{abstract}
This article is the second in the series and is devoted to the type $G_2$. The work consists of two parts.
In the first part we calculate the structure constants of the complex simple Lie algebra of type $G_2.$ 
All structure constants are represented as functions of the structure constants corresponding to extraspecial pairs.
The results obtained are used to calculate the commutator Chevalley formulas for $[x_r(u),x_s(y)],$ when the sum $r+s$ is a root.

Further, in the second part there is a table of structure constants and Chevalley commutator formulas in the special case,
when all structure constants corresponding to extraspecial pairs are positive.


\end{abstract}

\maketitle
\tableofcontents

\section{General case}

\subsection{Information about the Root System of Type $G_2$}

Let us choose in the three-dimensional Euclidean space with an orthonormal basis 
$\varepsilon_1,$ $\varepsilon_2,$ $\varepsilon_3$, as in \cite[p. 319]{Bur72}, two vectors
$$
a=\varepsilon_1-\varepsilon_2,\quad\ b=-2\varepsilon_1+\varepsilon_2+\varepsilon_3.
$$
Then the set of vectors
$$
-3a-2b,\ -3a-b,\ -2a-b,\ -a-b,\ -a,\ -b,\ a,\ b,\ a+b,\ 2a+b,\ 3a+b, \ 3a+2b
$$
forms the root system  of type $G_2,$ the roots
$$
a,\ b,\ a+b,\ 2a+b,\ 3a+b,\ 3a+2b
$$
form the subsystem of positive roots $G_2^+,$ and the roots $a$ and $b$ form
fundamental root system $\Pi(G_2).$
Ordering of positive roots $a\prec b\prec a+b\prec 2a+b\prec 3a+b\prec\ 3a+2b$
with respect to the linear order $\prec,$ determining the choice of positive roots.
Extraspecial pairs of roots are following:
$$
(a,b),\ (a,a+b),\ (a,2a+b),\ (b,3a+b).
$$

\subsection{Calculation of Structure Constants}

The following theorem calculates the structure constants of the complex simple Lie algebra of type $G_2.$
\medskip

\textbf{Theorem 1.} \textit{Let
$$
N_{a,b}=\epsilon_1,\ N_{a,a+b}=2\epsilon_2,\
N_{a,2a+b}=3\epsilon_3,\ N_{b,3a+b}=\epsilon_4
$$
and  $\epsilon_5=\frac{\epsilon_1\epsilon_3}{\epsilon_4}$.
Then the nonzero constants $N_{r,s}$ of the complex simple Lie algebra of type
$G_2$ have the form specified in following table:}

{\scriptsize
\begin{center}
\begin{tabular}{|r|r|r|r|r|r|r|r|r|r|r|r|r|}
\hline
 $N_{r,s}$ & $-3a-2b$      & $-3a-b$       & $-2a-b$        & $-a-b$         & $-a$           & $-b$          & $a$            & $b$           &\! $a+b$\!      &\! $2a+b$\!     &\! $3a+b$\!     &\! $3a+2b$\!   \\
\hline
 $-3a-2b$  &               &               &                &                &                &               &                & $\epsilon_4$  & $-\epsilon_5$  &  $\epsilon_5$  & $-\epsilon_4$  &               \\
\hline
 $-3a-b$   &               &               &                &                &                & $\epsilon_4$  & $\epsilon_3$   &               &                & $-\epsilon_3$  &                & $-\epsilon_4$ \\
\hline
 $-2a-b$   &               &               &                & $-3\epsilon_5$ & $3\epsilon_3$  &               & $2\epsilon_2$  &               & $-2\epsilon_2$ &                & $-\epsilon_3$ & $\epsilon_5$  \\
\hline
 $-a-b$    &               &               & $3\epsilon_5$  &                & $2\epsilon_2$  &               & $3\epsilon_1$  & $-\epsilon_1$ &                & $-2\epsilon_2$ &                & $-\epsilon_5$ \\
\hline
 $-b$      &               & $-\epsilon_4$ &                &                & $\epsilon_1$   &               &                &               & $-\epsilon_1$  &                &                & $\epsilon_4$  \\
\hline
 $-a$      &               &               & $-3\epsilon_3$ & $-2\epsilon_2$ &                & $-\epsilon_1$ &                &               & $3\epsilon_1$  & $2\epsilon_2$  & $\epsilon_3$   &   \\
\hline
 $a$       &               & $-\epsilon_3$ & $-2\epsilon_2$ & $-3\epsilon_1$ &                &               &                & $\epsilon_1$  & $2\epsilon_2$  & $3\epsilon_3$  &                &   \\
\hline
 $b$       & $-\epsilon_4$ &               &                & $\epsilon_1$   &                &               & $-\epsilon_1$  &               &                &                & $\epsilon_4$   &   \\
\hline
 $a+b$     & $\epsilon_5$  &               & $2\epsilon_2$  &                & $-3\epsilon_1$ & $\epsilon_1$  & $-2\epsilon_2$ &               &                & $-3\epsilon_5$ &                &   \\
\hline
 $2a+b$    & $-\epsilon_5$ & $\epsilon_3$  &                & $2\epsilon_2$  & $-2\epsilon_2$ &               & $-3\epsilon_3$ &               & $3\epsilon_5$  &                &                &   \\
\hline
 $3a+b$    & $\epsilon_4$  &               & $\epsilon_3$   &                & $-\epsilon_3$  &               &                & $-\epsilon_4$ &                &                &                &   \\
\hline
 $3a+2b$   &               & $\epsilon_4$  & $-\epsilon_5$  & $\epsilon_5$   &                & $-\epsilon_4$ &                &               &                &                &                &   \\
\hline
\end{tabular}
\end{center}
}
\medskip

\centerline{\textbf{Table 1.}}
\bigskip

\textbf{Proof.} According to \cite[Theorem 4.1.2]{Car72} the structure constants of a simple Lie algebra of type $\Phi$ over
$\mathbb{C}$ satisfy the following relations:

(i) $\displaystyle{N_{s,r}=-N_{r,s},\ r,s\in\Phi;}$
\bigskip

(ii) $\displaystyle{\frac{N_{r_1,r_2}}{(r_3,r_3)}=
\frac{N_{r_2,r_3}}{(r_1,r_1)}=\frac{N_{r_3,r_1}}{(r_2,r_2)},}$
\medskip

\noindent
if $\displaystyle{r_1,r_2,r_3\in\Phi}$ satisfy
$\displaystyle{r_1+r_2+r_3=0;}$
\bigskip

(iii) $\displaystyle{N_{r,s}N_{-r,-s}=-(p+1)^2,}$
\medskip

\noindent
if $r,s,r+s\in\Phi;$
\bigskip

(iv)
$\displaystyle{\frac{N_{r_1,r_2}N_{r_3,r_4}}{(r_1+r_2,r_1+r_2)}+
\frac{N_{r_2,r_3}N_{r_1,r_4}}{(r_2+r_3,r_2+r_3)}+
\frac{N_{r_3,r_1}N_{r_2,r_4}}{(r_3+r_1,r_3+r_1)}=0,}$
\medskip

\noindent
if $\displaystyle{r_1,r_2,r_3,r_4\in \Phi}$ satisfy
$\displaystyle{r_1+r_2+r_3+r_4=0}$ and if no pair opposite.
\medskip

From the equalities
$$
b+(-a-b)+a=0,\ (a+b)+(-2a-b)+a=0,\ (2a+b)+(-3a-b)+a=0,
$$
$$
(3a+b)+(-3a-2b)+b=0
$$
it follows that
$$
\frac{N_{b,-a-b}}{2}=\frac{N_{-a-b,a}}{6}=\frac{N_{a,b}}{2}=\frac{\epsilon_1} {2},
$$
$$
\frac{N_{a+b,-2a-b}}{2}=\frac{N_{-2a-b,a}}{2}=\frac{N_{a,a+b}}{2}=\epsilon_2,
$$
$$
\frac{N_{2a+b,-3a-b}}{2}=\frac{N_{-3a-b,a}}{2}=\frac{N_{a,2a+b}}{6}=\frac{\epsilon_3}{2},
$$
$$
\frac{N_{3a+b,-3a-2b}}{6}=\frac{N_{-3a-2b,b}}{6}=\frac{N_{b,3a+b}}{6}=\frac{\epsilon_4}{6}.
$$
Further, from the equality
$$
(a+b)+(2a+b)+(-b)+(-3a-b)=0
$$
it follows that
$$
\frac{N_{a+b,2a+b}N_{-b,-3a-b}}{6}+\frac{N_{-b,a+b}N_{2a+b,-3a-b }}{2}=0,
$$
where
$$
N_{a+b,2a+b}=-3\frac{N_{-b,a+b}N_{2a+b,-3a-b}}{N_{-b,-3a-b}}= -3\frac{(-\epsilon_1)\epsilon_3}{-\epsilon_4}
=-3\frac{\epsilon_1\epsilon_3}{\epsilon_4}=-3\epsilon_5.
$$
Finally, from the equality
$$
(2a+b)+(-3a-2b)+(a+b)=0
$$
it follows that
$$
\frac{N_{2a+b,-3a-2b}}{2}=\frac{N_{-3a-2b,a+b}}{2}=\frac{N_{a+b,2a+b }}{6}=-\frac{\epsilon_5}{2}.
$$

The remaining structure constants $N_{r,s}$ can be found from the relations
$$
N_{r,s}=N_{-s,-r}=-N_{s,r}=-N_{-r,-s}.
$$
The theorem is proved. \hfill $\blacksquare$
\medskip

\subsection{Chevalley's Commutator Formula}

Let $\Phi$ be a reduced indecomposable root system, $K$ is a
field (or associatively-commutative ring with identity). According to \cite[Theorem 5.2.2]{Car72} the commutator
$$
[x_s(u),x_r(t)]=x_s(u)^{-1}x_r(t)^{-1}x_s(u)x_r(t),
$$
where $r,s\in\Phi$ and $u,t\in K,$ is equal to one if $r+s\notin\Phi$ and $r\ne-s,$ and expands
into the product of root elements according to the formula
$$
\left[x_s(u),x_r(t)\right]=
 \prod_{{\scriptsize \begin{array}{c} ir+js\in\Phi,\\ i,j>0\end{array}}}x_{ir+js}(C_{ij,rs}( -t)^iu^j)
$$
if $r+s\in\Phi.$ The product is taken over all pairs of positive integers $i,j$ for which $ir+js$ is a root, in order of increasing $i+j.$ The constants $C_{ij,rs}$
are integers and are determined by the formulas \cite[theorem 5.2.2]{Car72}:
\begin{align*}
    C_{i1,rs} &= M_{r,s,i}, \\[3mm]
    C_{1j,rs} &= (-1)^j M_{s,r,j}, \\[3mm]
    C_{32,rs} &= \frac{1}{3}M_{r+s,r,2}, \\[3mm]
    C_{23,rs} &= -\frac{2}{3}M_{s+r,s,2}.
\end{align*}
In turn, the numbers $M_{r,s,i}$ are expressed in terms of structure
constants $N_{r,s}$ of the corresponding Lie algebra according to the formula \cite[p. 61]{Car72}
$$
M_{r,s,i}=\frac{1}{i!}N_{r,s}N_{r,r+s}\ldots N_{r,(i-1)r+s}.
$$
\vskip15pt

\centerline{\textbf{Positive roots}} \vskip10pt

\textit{Formula 1.}
$$
[x_b(u),x_a(t)]=
$$
$$
 x_{a+b}(C_{11,a,b}(-t)u)\, x_{2a+b}(C_{21,a,b}(-t)^2u)\, x_{3a+b}(C_{31,a,b}(-t)^3u)\,x_{3a+2b}(C_{32,a,b}(-t)^3u^2).
$$
Calculation of constants: \vskip3pt

$\displaystyle{C_{11,a,b}=M_{a,b,1}=N_{a,b}=\epsilon_1},$
\vskip7pt

$\displaystyle{C_{21,a,b}=M_{a,b,2}=\frac{1}{2!}N_{a,b}N_{a,a+b}=\frac{1}{2}(\epsilon_1)(2\epsilon_2)=\epsilon_1\epsilon_2
,}$ \vskip7pt

$\displaystyle{C_{31,a,b}=M_{a,b,3}=\frac{1}{3!}N_{a,b}N_{a,a+b}N_{a,2a+b}=\frac{1}{6}\epsilon_1(2\epsilon_2)(3\epsilon_3)=\epsilon_1\epsilon_2\epsilon_3
,}$ \vskip7pt

$\displaystyle{C_{32,a,b}=\frac{1}{3}M_{a+b,a,2}=\frac{1}{3!}N_{a+b,a}N_{a+b,2a+b}=\frac{1}{6}(-2\epsilon_2)(-3\epsilon_5)=\epsilon_2\epsilon_5
.}$ \vskip7pt

\noindent \textbf{Final formula:}
\vskip10pt

\begin{center}
 \fbox{$\displaystyle{[x_b(u),x_a(t)]=
 x_{a+b}(-\epsilon_1tu)\,x_{2a+b}(\epsilon_1\epsilon_2 t^2u)\,x_{3a+b}(-\epsilon_1\epsilon_2\epsilon_3 t^3u)\,x_{3a+2b}(-\epsilon_2\epsilon_5t^3u^2).}$}
\end{center}
\bigskip

\noindent\textit{Formula 2.}
$$
[x_a(u),x_b(t)]=
$$
$$
 x_{a+b}(C_{11,b,a}(-t)u)\, x_{2a+b}(C_{12,b,a}(-t)u^2)\, x_{3a+b}(C_{13,b,a}(-t)u^3)\,x_{3a+2b}(C_{23,b,a}(-t)^2u^3).
$$
Calculation of constants: \vskip3pt

$\displaystyle{C_{11,b,a}=M_{b,a,1}=N_{b,a}=-\epsilon_1},$
\vskip7pt

$\displaystyle{C_{12,b,a}=(-1)^2M_{a,b,2}=\frac{1}{2!}N_{a,b}N_{a,a+b}=\frac{1}{2}(\epsilon_1)(2\epsilon_2)=\epsilon_1\epsilon_2
,}$ \vskip7pt

$\displaystyle{C_{13,b,a}=(-1)^3M_{a,b,3}=-\frac{1}{3!}N_{a,b}N_{a,a+b}N_{a,2a+b}=-\frac{1}{6}\epsilon_1(2\epsilon_2)(3\epsilon_3)=-\epsilon_1\epsilon_2\epsilon_3
,}$ \vskip7pt

 $\displaystyle{C_{23,b,a}=-\frac23 M_{a+b,a,2}=-\frac{1}{3}N_{a+b,a}N_{a+b,2a+b}=-\frac{1}{3}(-2\epsilon_2)(-3\epsilon_5)=-2\epsilon_2\epsilon_5
.}$ \vskip7pt

\noindent \textbf{Final formula:}
\vskip10pt

\begin{center}
 \fbox{$\displaystyle{[x_a(u),x_b(t)]=
 x_{a+b}(\epsilon_1tu)\,x_{2a+b}(-\epsilon_1\epsilon_2 tu^2)\,x_{3a+b}(\epsilon_1\epsilon_2\epsilon_3 tu^3)\,x_{3a+2b}(-2\epsilon_2\epsilon_5t^2u^3).}$}
\end{center}
\bigskip

\noindent\textit{Formula 3.}
$$
[x_{a+b}(u),x_a(t)]=x_{2a+b}(C_{11,a,a+b}(-t)u)\,x_{3a+b}(C_{21,a,a+b}(-t)^2u)\,x_{3a+2b}(C_{12,a,a+b}(-t)u^2).
$$
Calculation of constants: \vskip3pt

$\displaystyle{C_{11,a,a+b}=M_{a,a+b,1}=N_{a,a+b}=2\epsilon_2},$
\vskip7pt

$\displaystyle{C_{21,a,a+b}=M_{a,a+b,2}=\frac{1}{2!}N_{a,a+b}N_{a,2a+b}=\frac{1}{2}(2\epsilon_2)(3\epsilon_3)=3\epsilon_2\epsilon_3
,}$ \vskip7pt

 $\displaystyle{C_{12,a,a+b}=(-1)^2 M_{a+b,a,2}=\frac{1}{2!}N_{a+b,a}N_{a+b,2a+b}=\frac{1}{2}(-2\epsilon_2)(-3\epsilon_5)=3\epsilon_2\epsilon_5.
}$ \vskip7pt

\noindent \textbf{Final formula:}
\vskip10pt

\begin{center}
 \fbox{$\displaystyle{[x_{a+b}(u),x_a(t)]=x_{2a+b}(-2\epsilon_2tu)\,x_{3a+b}(3\epsilon_2\epsilon_3t^2u)\,x_{3a+2b}(-3\epsilon_2\epsilon_5tu^2).}$}
\end{center}
\bigskip

\noindent\textit{Formula 4.}
$$
[x_{2a+b}(u),x_a(t)]=x_{3a+b}(C_{11,a,2a+b}(-t)u).
$$
Calculation of constants: \vskip3pt

$\displaystyle{C_{11,a,2a+b}=M_{a,2a+b,1}=N_{a,2a+b}=3\epsilon_3}.$
\vskip7pt

\noindent \textbf{Final formula:}
\vskip10pt

\begin{center}
 \fbox{$\displaystyle{[x_{2a+b}(u),x_a(t)]=x_{3a+b}(-3\epsilon_3tu).}$}
\end{center}
\bigskip

\noindent\textit{Formula 5.}
$$
[x_{3a+b}(u),x_b(t)]=x_{3a+2b}(C_{11,b,3a+b}(-t)u).
$$
Calculation of constants: \vskip3pt

$\displaystyle{C_{11,b,3a+b}=M_{b,3a+b,1}=N_{b,3a+b}=\epsilon_4}.$
\vskip7pt

\noindent \textbf{Final formula:}
\vskip10pt

\begin{center}
 \fbox{$\displaystyle{[x_{3a+b}(u),x_b(t)]=x_{3a+2b}(-\epsilon_4tu).}$}
\end{center}
\bigskip

\noindent\textit{Formula 6.}
$$
[x_{2a+b}(u),x_{a+b}(t)]=x_{3a+2b}(C_{11,a+b,2a+b}(-t)u).
$$
Calculation of constants: \vskip3pt

$\displaystyle{C_{11,a+b,2a+b}=M_{a+b,2a+b,1}=N_{a+b,2a+b}=-3\epsilon_5}.$
\vskip7pt

\noindent \textbf{Final formula:}
\vskip10pt

\begin{center}
 \fbox{$\displaystyle{[x_{2a+b}(u),x_{a+b}(t)]=x_{3a+2b}(3\epsilon_5tu).}$}
\end{center}

\vskip15pt \centerline{\textbf{Negative roots}} \vskip15pt

\noindent\textit{Formula 7.}
$$
[x_{-b}(u),x_{-a}(t)]=x_{-a-b}(C_{11,-a,-b}(-t)u)\times
$$
$$
x_{-2a-b}(C_{21,-a,-b}(-t)^2u)\, x_{-3a-b}(C_{31,-a,-b}(-t)^3u)\,x_{-3a-2b}(C_{32,-a,-b}(-t)^3u^2).
$$
Calculation of constants: \vskip3pt

$\displaystyle{C_{11,-a,-b}=M_{-a,-b,1}=N_{-a,-b}=-\epsilon_1},$
\vskip7pt

$\displaystyle{C_{21,-a,-b}=M_{-a,-b,2}=\frac{1}{2!}N_{-a,-b}N_{-a,-a-b}=\frac{1}{2}(-\epsilon_1)(-2\epsilon_2)=\epsilon_1\epsilon_2
,}$ \vskip7pt

$\displaystyle{C_{31,-a,-b}=M_{-a,-b,3}=\frac{1}{3!}N_{-a,-b}N_{-a,-a-b}N_{-a,-2a-b}=\frac{1}{6}(-\epsilon_1)(-2\epsilon_2)(-3\epsilon_3)=-\epsilon_1\epsilon_2\epsilon_3
,}$ \vskip7pt

$\displaystyle{C_{32,-a,-b}=\frac{1}{3}M_{-a-b,-a,2}=\frac{1}{3!}N_{-a-b,-a}N_{-a-b,-2a-b}=\frac{1}{6}(2\epsilon_2)(3\epsilon_5)=\epsilon_2\epsilon_5
.}$ \vskip7pt

\noindent \textbf{Final formula:}
\vskip10pt

\begin{center}
 \fbox{$\displaystyle{[x_{-b}(u),x_{-a}(t)]=
 x_{-a-b}(\epsilon_1tu)\,x_{-2a-b}(\epsilon_1\epsilon_2 t^2u)\,x_{-3a-b}(\epsilon_1\epsilon_2\epsilon_3 t^3u)\,x_{-3a-2b}(-\epsilon_2\epsilon_5t^3u^2).}$}
\end{center}
\bigskip

\noindent\textit{Formula 8.}
$$
[x_{-a}(u),x_{-b}(t)]=x_{-a-b}(C_{11,-b,-a}(-t)u)\times
$$
$$
x_{-2a-b}(C_{12,-b,-a}(-t)u^2)\, x_{-3a-b}(C_{13,-b,-a}(-t)u^3)\,x_{-3a-2b}(C_{23,-b,-a}(-t)^2u^3).
$$
Calculation of constants: \vskip3pt

 $\displaystyle{C_{11,-b,-a}=M_{-b,-a,1}=N_{-b,-a}=\epsilon_1},$
\vskip7pt

 $\displaystyle{C_{12,-b,-a}=(-1)^2M_{-a,-b,2}=\frac{1}{2!}N_{-a,-b}N_{-a,-a-b}=\frac{1}{2}(-\epsilon_1)(-2\epsilon_2)=\epsilon_1\epsilon_2
,}$ \vskip7pt

 $\displaystyle{C_{13,-b,-a}=(-1)^3M_{-a,-b,3}=-\frac{1}{3!}N_{-a,-b}N_{-a,-a-b}N_{-a,-2a-b}=-\frac{1}{6}(-\epsilon_1)(-2\epsilon_2)(-3\epsilon_3)=}$
\vskip7pt

$\displaystyle{\ \ \ \ \ \ \ \ \ \ \ \ =\epsilon_1\epsilon_2\epsilon_3,}$ 
\vskip7pt

 $\displaystyle{C_{23,-b,-a}=-\frac23 M_{-a-b,-a,2}=-\frac{1}{3}N_{-a-b,-a}N_{-a-b,-2a-b}=-\frac{1}{3}(2\epsilon_2)(3\epsilon_5)=-2\epsilon_2\epsilon_5
.}$ \vskip7pt

\noindent \textbf{Final formula:}
\vskip10pt

\begin{center}
 \fbox{$\displaystyle{[x_{-a}(u),x_{-b}(t)]=
 x_{-a-b}(-\epsilon_1tu)\,x_{-2a-b}(-\epsilon_1\epsilon_2 tu^2)x_{-3a-b}(-\epsilon_1\epsilon_2\epsilon_3 tu^3)
 x_{-3a-2b}(-2\epsilon_2\epsilon_5t^2u^3).}$}
\end{center}
\bigskip

\noindent\textit{Formula 9.}
$$
[x_{-a-b}(u),x_{-a}(t)]=x_{-2a-b}(C_{11,-a,-a-b}(-t)u)\times
$$
$$
x_{-3a-b}(C_{21,-a,-a-b}(-t)^2u)\,x_{-3a-2b}(C_{12,-a,-a-b}(-t)u^2).
$$
Calculation of constants: \vskip3pt

$\displaystyle{C_{11,-a,-a-b}=M_{-a,-a-b,1}=N_{-a,-a-b}=-2\epsilon_2},$
\vskip7pt

$\displaystyle{C_{21,-a,-a-b}=M_{-a,-a-b,2}=\frac{1}{2!}N_{-a,-a-b}N_{-a,-2a-b}=\frac{1}{2}(-2\epsilon_2)(-3\epsilon_3)=3\epsilon_2\epsilon_3
,}$ \vskip7pt

 $\displaystyle{C_{12,-a,-a-b}=(-1)^2 M_{-a-b,-a,2}=\frac{1}{2!}N_{-a-b,-a}N_{-a-b,-2a-b}=\frac{1}{2}(2\epsilon_2)(3\epsilon_5)=3\epsilon_2\epsilon_5.
}$ \vskip7pt

\noindent \textbf{Final formula:}
\vskip10pt

\begin{center}
 \fbox{$\displaystyle{[x_{-a-b}(u),x_{-a}(t)]=x_{-2a-b}(2\epsilon_2tu)\,x_{-3a-b}(3\epsilon_2\epsilon_3t^2u)\,x_{-3a-2b}(-3\epsilon_2\epsilon_5tu^2).}$}
\end{center}
\bigskip

\noindent\textit{Formula 10.}
$$
[x_{-2a-b}(u),x_{-a}(t)]=x_{-3a-b}(C_{11,-a,-2a-b}(-t)u).
$$
Calculation of constants: \vskip3pt

$\displaystyle{C_{11,-a,-2a-b}=M_{-a,-2a-b,1}=N_{-a,-2a-b}=-3\epsilon_3}.$
\vskip7pt

\noindent \textbf{Final formula:}
\vskip10pt

\begin{center}
 \fbox{$\displaystyle{[x_{-2a-b}(u),x_{-a}(t)]=x_{-3a-b}(3\epsilon_3tu).}$}
\end{center}
\bigskip

\noindent\textit{Formula 11.}
$$
[x_{-3a-b}(u),x_{-b}(t)]=x_{-3a-2b}(C_{11,-b,-3a-b}(-t)u).
$$
Calculation of constants: \vskip3pt

$\displaystyle{C_{11,-b,-3a-b}=M_{-b,-3a-b,1}=N_{-b,-3a-b}=-\epsilon_4}.$
\vskip7pt

\noindent \textbf{Final formula:}
\vskip10pt

\begin{center}
 \fbox{$\displaystyle{[x_{-3a-b}(u),x_{-b}(t)]=x_{-3a-2b}(\epsilon_4tu).}$}
\end{center}
\bigskip

\noindent\textit{Formula 12.}
$$
[x_{-2a-b}(u),x_{-a-b}(t)]=x_{-3a-2b}(C_{11,-a-b,-2a-b}(-t)u).
$$
Calculation of constants: \vskip3pt

$\displaystyle{C_{11,-a-b,-2a-b}=M_{-a-b,-2a-b,1}=N_{-a-b,-2a-b}=3\epsilon_5}.$
\vskip7pt

\noindent \textbf{Final formula:}
\vskip10pt

\begin{center}
 \fbox{$\displaystyle{[x_{-2a-b}(u),x_{-a-b}(t)]=x_{-3a-2b}(-3\epsilon_5tu).}$}
\end{center}

\vskip15pt \centerline{\textbf{Positive and negative roots}} \vskip15pt

\noindent\textit{Formula 13.}
$$
[x_{-a-b}(u),x_{a}(t)]=x_{-b}(C_{11,a,-a-b}(-t)u).
$$
Calculation of constants: \vskip3pt

$\displaystyle{C_{11,a,-a-b}=M_{a,-a-b,1}=N_{a,-a-b}=-3\epsilon_1}.$
\vskip7pt

\noindent \textbf{Final formula:}
\vskip10pt

\begin{center}
 \fbox{$\displaystyle{[x_{-a-b}(u),x_{a}(t)]=x_{-b}(3\epsilon_1tu).}$}
\end{center}
\bigskip

\noindent\textit{Formula 14.}
$$
[x_{-2a-b}(u),x_a(t)]=
$$
$$
=x_{-a-b}(C_{11,a,-2a-b}(-t)u)\,x_{-b}(C_{21,a,-2a-b}(-t)^2u)\,x_{-3a-2b}(C_{12,a,-2a-b}(-t)u^2).
$$
Calculation of constants: \vskip3pt

$\displaystyle{C_{11,a,-2a-b}=M_{a,-2a-b,1}=N_{a,-2a-b}=-2\epsilon_2},$
\vskip7pt

$\displaystyle{C_{21,a,-2a-b}=M_{a,-2a-b,2}=\frac{1}{2!}N_{a,-2a-b}N_{a,-a-b}=\frac{1}{2}(-2\epsilon_2)(-3\epsilon_1)=3\epsilon_1\epsilon_2
,}$ \vskip7pt

 $\displaystyle{C_{12,a,-2a-b}=(-1)^2 M_{-2a-b,a,2}=\frac{1}{2!}N_{-2a-b,a}N_{-2a-b,-a-b}=\frac{1}{2}(2\epsilon_2)(-3\epsilon_5)=-3\epsilon_2\epsilon_5.
}$ \vskip7pt

\noindent \textbf{Final formula:}
\vskip10pt

\begin{center}
 \fbox{$\displaystyle{[x_{-2a-b}(u),x_a(t)]=x_{-a-b}(2\epsilon_2tu)\,x_{-b}(3\epsilon_1\epsilon_2t^2u)\,x_{-3a-2b}(3\epsilon_2\epsilon_5tu^2).}$}
\end{center}
\bigskip

\noindent\textit{Formula 15.}
$$
[x_{-3a-b}(u),x_a(t)]=x_{-2a-b}(C_{11,a,-3a-b}(-t)u)\times
$$
$$
x_{-a-b}(C_{21,a,-3a-b}(-t)^2u)\, x_{-b}(C_{31,a,-3a-b}(-t)^3u)\,x_{-3a-2b}(C_{32,a,-3a-b}(-t)^3u^2).
$$
Calculation of constants: \vskip3pt

$\displaystyle{C_{11,a,-3a-b}=M_{a,-3a-b,1}=N_{a,-3a-b}=-\epsilon_3},$
\vskip7pt

$\displaystyle{C_{21,a,-3a-b}=M_{a,-3a-b,2}=\frac{1}{2!}N_{a,-3a-b}N_{a,-2a-b}=\frac{1}{2}(-\epsilon_3)(-2\epsilon_2)=\epsilon_2\epsilon_3
,}$ \vskip7pt

$\displaystyle{C_{31,a,-3a-b}=M_{a,-3a-b,3}=\frac{1}{3!}N_{a,-3a-b}N_{a,-2a-b}N_{a,-a-b}=\frac{1}{6}(-\epsilon_3)(-2\epsilon_2)(-3\epsilon_1)=}$
\vskip7pt

$\displaystyle{\ \ \ \ \ \ \ \ \ \ \ \ \ \ \, \! =-\epsilon_1\epsilon_2\epsilon_3,}$ 
\vskip7pt

$\displaystyle{C_{32,a,-3a-b}=\frac{1}{3}M_{-2a-b,a,2}=\frac{1}{3!}N_{-2a-b,a}N_{-2a-b,-a-b}=
                              \frac{1}{6}(2\epsilon_2)(-3\epsilon_5)=-\epsilon_2\epsilon_5.}$
\vskip7pt

\noindent \textbf{Final formula:}
\vskip10pt

\begin{center}
 \fbox{$\displaystyle{[x_{-3a-b}(u),x_a(t)]=
 x_{-2a-b}(\epsilon_3tu)\,x_{-a-b}(\epsilon_2\epsilon_3 t^2u)\,x_{-b}(\epsilon_1\epsilon_2\epsilon_3 t^3u)\,x_{-3a-2b}(\epsilon_2\epsilon_5t^3u^2).}$}
\end{center}
\bigskip

\noindent\textit{Formula 16.}
$$
[x_a(u),x_{-3a-b}(t)]=x_{-2a-b}(C_{11,-3a-b,a}(-t)u)\times
$$
$$
 x_{-a-b}(C_{12,-3a-b,a}(-t)u^2)\, x_{-b}(C_{13,-3a-b,a}(-t)u^3)\,x_{-3a-2b}(C_{23,-3a-b,a}(-t)^2u^3).
$$
Calculation of constants: \vskip3pt

$\displaystyle{C_{11,-3a-b,a}=M_{-3a-b,a,1}=N_{-3a-b,a}=\epsilon_3},$
\vskip7pt

$\displaystyle{C_{12,-3a-b,a}=(-1)^2M_{a,-3a-b,2}=\frac{1}{2!}N_{a,-3a-b}N_{a,-2a-b}=\frac{1}{2}(-\epsilon_3)(-2\epsilon_2)=\epsilon_2\epsilon_3
,}$ \vskip7pt

$\displaystyle{C_{31,a,-3a-b}=(-1)^3M_{a,-3a-b,3}=-\frac{1}{3!}N_{a,-3a-b}N_{a,-2a-b}N_{a,-a-b}=}$

$\qquad\qquad\ \ \displaystyle{=-\frac{1}{6}(-\epsilon_3)(-2\epsilon_2)(-3\epsilon_1)=\epsilon_1\epsilon_2\epsilon_3,}$
\vskip7pt

$\displaystyle{C_{23,-3a-b,a}=-\frac{2}{3}M_{-2a-b,a,2}=-\frac{1}{3}N_{-2a-b,a}N_{-2a-b,-a-b}=-\frac{1}{3}(2\epsilon_2)(-3\epsilon_5)=
2\epsilon_2\epsilon_5.}$ \vskip7pt

\noindent \textbf{Final formula:}
\vskip10pt

\begin{center}
 \fbox{$\displaystyle{[x_{a}(u),x_{-3a-b}(t)]=
 x_{-2a-b}(-\epsilon_3tu)\,x_{-a-b}(-\epsilon_2\epsilon_3 tu^2)\,x_{-b}(-\epsilon_1\epsilon_2\epsilon_3 tu^3)\,x_{-3a-2b}(2\epsilon_2\epsilon_5t^2u^3).}$}
\end{center}
\bigskip

\noindent\textit{Formula 17.}
$$
[x_b(u),x_{-a-b}(t)]=x_{-a}(C_{11,-a-b,b}(-t)u)\times
$$
$$
x_{-2a-b}(C_{21,-a-b,b}(-t)^2u)\,x_{-3a-2b}(C_{31,-a-b,b}(-t)^3u)\,x_{-3a-b}(C_{32,-a-b,b}(-t)^3u^2).
$$
Calculation of constants: \vskip3pt

$\displaystyle{C_{11,-a-b,b}=M_{-a-b,b,1}=N_{-a-b,b}=-\epsilon_1},$
\vskip7pt

$\displaystyle{C_{21,-a-b,b}=M_{-a-b,b,2}=\frac{1}{2!}N_{-a-b,b}N_{-a-b,-a}=\frac{1}{2}(-\epsilon_1)(2\epsilon_2)=-\epsilon_1\epsilon_2
,}$ \vskip7pt

$\displaystyle{C_{31,-a-b,b}=M_{-a-b,b,3}=\frac{1}{3!}N_{-a-b,b}N_{-a-b,-a}N_{-a-b,-2a-b}=\frac{1}{6}(-\epsilon_1)(2\epsilon_2)(3\epsilon_5)=
-\epsilon_1\epsilon_2\epsilon_5,}$ \vskip7pt

$\displaystyle{C_{32,-a-b,b}=\frac{1}{3}M_{-a,-a-b,2}=\frac{1}{3!}N_{-a,-a-b}N_{-a,-2a-b}=\frac{1}{6}(-2\epsilon_2)(-3\epsilon_3)=
\epsilon_2\epsilon_3.}$ \vskip7pt

\noindent \textbf{Final formula:}
\vskip10pt

\begin{center}
 \fbox{$\displaystyle{[x_b(u),x_{-a-b}(t)]=x_{-a}(\epsilon_1tu)\,x_{-2a-b}(-\epsilon_1\epsilon_2t^2u)\,x_{-3a-2b}(\epsilon_1\epsilon_2\epsilon_5t^3u)\,
x_{-3a-b}(-\epsilon_2\epsilon_3t^3u^2).}$}
\end{center}
\bigskip

\noindent\textit{Formula 18.}
$$
[x_{-a-b}(u),x_b(t)]=x_{-a}(C_{11,b,-a-b}(-t)u)\times
$$
$$
x_{-2a-b}(C_{12,b,-a-b}(-t)u^2)\,x_{-3a-2b}(C_{13,b,-a-b}(-t)u^3)\,x_{-3a-b}(C_{23,b,-a-b}(-t)^2u^3).
$$
Calculation of constants: \vskip3pt

$\displaystyle{C_{11,b,-a-b}=M_{b,-a-b,1}=N_{b,-a-b}=\epsilon_1},$
\vskip7pt

$\displaystyle{C_{12,b,-a-b}=(-1)^2M_{-a-b,b,2}=\frac{1}{2!}N_{-a-b,b}N_{-a-b,-a}=\frac{1}{2}(-\epsilon_1)(2\epsilon_2)=-\epsilon_1\epsilon_2
,}$ \vskip7pt

$\displaystyle{C_{13,b,-a-b}=(-1)^3M_{-a-b,b,3}=-\frac{1}{3!}N_{-a-b,b}N_{-a-b,-a}N_{-a-b,-2a-b}=-\frac{1}{6}(-\epsilon_1)(2\epsilon_2)(3\epsilon_5)=}$
\vskip7pt

$\displaystyle{\ \ \ \ \ \ \ \ \ \ \ \ \ =\epsilon_1\epsilon_2\epsilon_5,}$ \vskip7pt

$\displaystyle{C_{23,b,-a-b}=-\frac{2}{3}M_{-a,-a-b,2}=-\frac{1}{3}N_{-a,-a-b}N_{-a,-2a-b}=-\frac{1}{3}(-2\epsilon_2)(-3\epsilon_3)=
-2\epsilon_2\epsilon_3.}$ \vskip7pt

\noindent \textbf{Final formula:}
\vskip10pt

\begin{center}
 \fbox{$\displaystyle{[x_{-a-b}(u),x_b(t)]=x_{-a}(-\epsilon_1tu)\,x_{-2a-b}(\epsilon_1\epsilon_2tu^2)\,x_{-3a-2b}(-\epsilon_1\epsilon_2\epsilon_5tu^3)\,
x_{-3a-b}(-2\epsilon_2\epsilon_3t^2u^3).}$}
\end{center}
\bigskip

\noindent\textit{Formula 19.}
$$
[x_{-3a-2b}(u),x_b(t)]=x_{-3a-b}(C_{11,b,-3a-2b}(-t)u).
$$

Calculation of constants: \vskip3pt

$\displaystyle{C_{11,b,-3a-2b}=M_{b,-3a-2b,1}=N_{b,-3a-2b}=-\epsilon_4}.$
\vskip7pt

\noindent \textbf{Final formula:}
\vskip10pt

\begin{center}
 \fbox{$\displaystyle{[x_{-3a-2b}(u),x_b(t)]=x_{-3a-b}(\epsilon_4tu).}$}
\end{center}
\bigskip

\noindent\textit{Formula 20.}
$$
[x_{-a}(u),x_{a+b}(t)]=x_b(C_{11,a+b,-a}(-t)u).
$$

Calculation of constants: \vskip3pt

$\displaystyle{C_{11,a+b,-a}=M_{a+b,-a,1}=N_{a+b,-a}=-3\epsilon_1}.$
\vskip7pt

\noindent \textbf{Final formula:}
\vskip10pt

\begin{center}
 \fbox{$\displaystyle{[x_{-a}(u),x_{a+b}(t)]=x_b(3\epsilon_1tu).}$}
\end{center}
\bigskip

\noindent\textit{Formula 21.}
$$
[x_{-b}(u),x_{a+b}(t)]=x_a(C_{11,a+b,-b}(-t)u)\times
$$
$$
x_{2a+b}(C_{21,a+b,-b}(-t)^2u)\,x_{3a+2b}(C_{31,a+b,-b}(-t)^3u)\,x_{3a+b}(C_{32,a+b,-b}(-t)^3u^2).
$$
Calculation of constants: \vskip3pt

$\displaystyle{C_{11,a+b,-b}=M_{a+b,-b,1}=N_{a+b,-b}=\epsilon_1},$
\vskip7pt

$\displaystyle{C_{21,a+b,-b}=M_{a+b,-b,2}=\frac{1}{2!}N_{a+b,-b}N_{a+b,a}=\frac{1}{2}(\epsilon_1)(-2\epsilon_2)=-\epsilon_1\epsilon_2
,}$ \vskip7pt

$\displaystyle{C_{31,a+b,-b}=M_{a+b,-b,3}=\frac{1}{3!}N_{a+b,-b}N_{a+b,a}N_{a+b,2a+b}=\frac{1}{6}(\epsilon_1)(-2\epsilon_2)(-3\epsilon_5)=
\epsilon_1\epsilon_2\epsilon_5,}$ \vskip7pt

$\displaystyle{C_{32,a+b,-b}=\frac{1}{3}M_{a,a+b,2}=\frac{1}{3!}N_{a,a+b}N_{a,2a+b}=\frac{1}{6}(2\epsilon_2)(3\epsilon_3)=
\epsilon_2\epsilon_3.}$ \vskip7pt

\noindent \textbf{Final formula:}
\vskip10pt

\begin{center}
 \fbox{$\displaystyle{[x_{-b}(u),x_{a+b}(t)]=x_a(-\epsilon_1tu)\,x_{2a+b}(-\epsilon_1\epsilon_2t^2u)\,x_{3a+2b}(-\epsilon_1\epsilon_2\epsilon_5t^3u)\,
x_{3a+b}(-\epsilon_2\epsilon_3t^3u^2).}$}
\end{center}
\bigskip

\noindent\textit{Formula 22.}
$$
[x_{a+b}(u),x_{-b}(t)]=x_a(C_{11,-b,a+b}(-t)u)\times
$$
$$
x_{2a+b}(C_{12,-b,a+b}(-t)u^2)\,x_{3a+2b}(C_{13,-b,a+b}(-t)u^3)\,x_{3a+b}(C_{23,-b,a+b}(-t)^2u^3).
$$
Calculation of constants: \vskip3pt

$\displaystyle{C_{11,-b,a+b}=M_{-b,a+b,1}=N_{-b,a+b}=-\epsilon_1},$
\vskip7pt

$\displaystyle{C_{12,-b,a+b}=(-1)^2M_{a+b,-b,2}=\frac{1}{2!}N_{a+b,-b}N_{a+b,a}=\frac{1}{2}(\epsilon_1)(-2\epsilon_2)=-\epsilon_1\epsilon_2
,}$ \vskip7pt

$\displaystyle{C_{13,-b,a+b}=(-1)^3M_{a+b,-b,3}=-\frac{1}{3!}N_{a+b,-b}N_{a+b,a}N_{a+b,2a+b}=-\frac{1}{6}\epsilon_1(-2\epsilon_2)(-3\epsilon_5)=}$
\vskip7pt

$\displaystyle{\ \ \ \ \ \ \ \ \ \ \ \ \ =-\epsilon_1\epsilon_2\epsilon_5,}$ 
\vskip7pt

$\displaystyle{C_{23,-b,a+b}=-\frac{2}{3}M_{a,a+b,2}=-\frac{1}{3}N_{a,a+b}N_{a,2a+b}=-\frac{1}{3}(2\epsilon_2)(3\epsilon_3)=
-2\epsilon_2\epsilon_3.}$ \vskip7pt

\noindent \textbf{Final formula:}
\vskip10pt

\begin{center}
 \fbox{$\displaystyle{[x_{a+b}(u),x_{-b}(t)]=x_a(\epsilon_1tu)\,x_{2a+b}(\epsilon_1\epsilon_2tu^2)\,x_{3a+2b}(\epsilon_1\epsilon_2\epsilon_5tu^3)\,
x_{3a+b}(-2\epsilon_2\epsilon_3t^2u^3).}$}
\end{center}
\bigskip

\noindent\textit{Formula 23.}
$$
[x_{-2a-b}(u),x_{a+b}(t)]=
$$
$$
x_{-a}(C_{11,a+b,-2a-b}(-t)u)\, x_b(C_{21,a+b,-2a-b}(-t)^2u)\,x_{-3a-b}(C_{12,a+b,-2a-b}(-t)u^2).
$$
Calculation of constants: \vskip3pt

$\displaystyle{C_{11,a+b,-2a-b}=M_{a+b,-2a-b,1}=N_{a+b,-2a-b}=2\epsilon_2},$
\vskip7pt

$\displaystyle{C_{21,a+b,-2a-b}=M_{a+b,-2a-b,2}=\frac{1}{2!}N_{a+b,-2a-b}N_{a+b,-a}=\frac{1}{2}(2\epsilon_2)(-3\epsilon_1)=-3\epsilon_1\epsilon_2
,}$ \vskip7pt

$\displaystyle{C_{12,a+b,-2a-b}=(-1)^2M_{-2a-b,a+b,2}=\frac{1}{2!}N_{-2a-b,a+b}N_{-2a-b,-a}=\frac{1}{2!}(-2\epsilon_2)(3\epsilon_3)=
-3\epsilon_2\epsilon_3.}$ \vskip7pt

\noindent \textbf{Final formula:}
\vskip10pt

\begin{center}
 \fbox{$\displaystyle{[x_{-2a-b}(u),x_{a+b}(t)]=x_{-a}(-2\epsilon_2tu)\,x_b(-3\epsilon_1\epsilon_2t^2u)\,x_{-3a-b}(3\epsilon_2\epsilon_3tu^2).}$}
\end{center}
\bigskip

\noindent\textit{Formula 24.}
$$
[x_{-3a-2b}(u),x_{a+b}(t)]=x_{-2a-b}(C_{11,a+b,-3a-2b}(-t)u)\times
$$
$$
x_{-a}(C_{21,a+b,-3a-2b}(-t)^2u)\,x_{b}(C_{31,a+b,-3a-2b}(-t)^3u)\,x_{-3a-b}(C_{32,a+b,-3a-2b}(-t)^3u^2).
$$
Calculation of constants: \vskip3pt

$\displaystyle{C_{11,a+b,-3a-2b}=M_{a+b,-3a-2b,1}=N_{a+b,-3a-2b}=\epsilon_5},$
\vskip7pt

$\displaystyle{C_{21,a+b,-3a-2b}=M_{a+b,-3a-2b,2}=\frac{1}{2!}N_{a+b,-3a-2b}N_{a+b,-2a-b}=\frac{1}{2}(\epsilon_5)(2\epsilon_2)=\epsilon_2\epsilon_5
,}$ \vskip7pt

$\displaystyle{C_{31,a+b,-3a-2b}=M_{a+b,-3a-2b,3}=\frac{1}{3!}N_{a+b,-3a-2b}N_{a+b,-2a-b}N_{a+b,-a}=\frac{1}{6}(\epsilon_5)(2\epsilon_2)(-3\epsilon_1)=}
$

$\displaystyle{\ \ \ \ \ \ \ \ \ \ \ \ \ \ \ \ \ \
=-\epsilon_1\epsilon_2\epsilon_5,}$ \vskip7pt

$\displaystyle{C_{32,a+b,-3a-2b}=\frac{1}{3}M_{-2a-b,a+b,2}=\frac{1}{3!}N_{-2a-b,a+b}N_{-2a-b,-a}=\frac{1}{6}(-2\epsilon_2)(3\epsilon_3)=
-\epsilon_2\epsilon_3.}$ \vskip7pt

\noindent \textbf{Final formula:}
\vskip10pt

\begin{center}
 \fbox{$\displaystyle{[x_{-3a-2b}(u),x_{a+b}(t)]=x_{-2a-b}(-\epsilon_5tu)\,x_{-a}(\epsilon_2\epsilon_5t^2u)\,x_b(\epsilon_1\epsilon_2\epsilon_5t^3u)\,
x_{-3a-b}(\epsilon_2\epsilon_3t^3u^2).}$}
\end{center}
\bigskip

\noindent\textit{Formula 25.}
$$
[x_{a+b}(u),x_{-3a-2b}(t)]=x_{-2a-b}(C_{11,-3a-2b,a+b}(-t)u)\times
$$
$$
x_{-a}(C_{12,-3a-2b,a+b}(-t)u^2)\,x_{b}(C_{13,-3a-2b,a+b}(-t)u^3)\,x_{-3a-b}(C_{23,-3a-2b,a+b}(-t)^2u^3).
$$
Calculation of constants: \vskip3pt

$\displaystyle{C_{11,-3a-2b,a+b}=M_{-3a-2b,a+b,1}=N_{-3a-2b,a+b}=-\epsilon_5},$
\vskip7pt

$\displaystyle{C_{12,-3a-2b,a+b}=(-1)^2M_{a+b,-3a-2b,2}=\frac{1}{2!}N_{a+b,-3a-2b}N_{a+b,-2a-b}=\frac{1}{2}(\epsilon_5)(2\epsilon_2)=\epsilon_2\epsilon_5
,}$ \vskip7pt

$\displaystyle{C_{13,-3a-2b,a+b}=(-1)^3M_{a+b,-3a-2b,3}=-\frac{1}{3!}N_{a+b,-3a-2b}N_{a+b,-2a-b}N_{a+b,-a}=}$\vskip7pt

$\displaystyle{\ \ \ \ \ \ \ \ \ \ \ \ \ \ \ \ \ \
=-\frac{1}{6}(\epsilon_5)(2\epsilon_2)(-3\epsilon_1)=\epsilon_1\epsilon_2\epsilon_5,}$
\vskip7pt

$\displaystyle{C_{23,-3a-2b,a+b}=-\frac{2}{3}M_{-2a-b,a+b,2}=-\frac{1}{3}N_{-2a-b,a+b}N_{-2a-b,-a}=-\frac{1}{3}(-2\epsilon_2)(3\epsilon_3)=
2\epsilon_2\epsilon_3.}$ \vskip7pt

\noindent \textbf{Final formula:}
\vskip10pt

\begin{center}
 \fbox{$\displaystyle{[x_{a+b}(u),x_{-3a-2b}(t)]=x_{-2a-b}(\epsilon_5tu)\,x_{-a}(-\epsilon_2\epsilon_5tu^2)\,x_b(-\epsilon_1\epsilon_2\epsilon_5tu^3)\,
x_{-3a-b}(2\epsilon_2\epsilon_3t^2u^3).}$}
\end{center}
\bigskip

\noindent\textit{Formula 26.}
$$
[x_{-a}(u),x_{2a+b}(t)]=x_{a+b}(C_{11,2a+b,-a}(-t)u)
 x_{3a+2b}(C_{21,2a+b,-a}(-t)^2u)\,x_b(C_{12,2a+b,-a}(-t)u^2).
$$
Calculation of constants: \vskip3pt

$\displaystyle{C_{11,2a+b,-a}=M_{2a+b,-a,1}=N_{2a+b,-a}=-2\epsilon_2},$
\vskip7pt

$\displaystyle{C_{21,2a+b,-a}=M_{2a+b,-a,2}=\frac{1}{2!}N_{2a+b,-a}N_{2a+b,a+b}=\frac{1}{2}(-2\epsilon_2)(3\epsilon_5)=-3\epsilon_2\epsilon_5
,}$ \vskip7pt

$\displaystyle{C_{12,2a+b,-a}=(-1)^2M_{-a,2a+b,2}=\frac{1}{2!}N_{-a,2a+b}N_{-a,a+b}=\frac{1}{2}(2\epsilon_2)(3\epsilon_1)=
3\epsilon_1\epsilon_2.}$ \vskip7pt

\noindent \textbf{Final formula:}
\vskip10pt

\begin{center}
 \fbox{$\displaystyle{[x_{-a}(u),x_{2a+b}(t)]=x_{a+b}(2\epsilon_2tu)\,x_{3a+2b}(-3\epsilon_2\epsilon_5t^2u)\,x_b(-3\epsilon_1\epsilon_2tu^2).}$}
\end{center}
\bigskip

\noindent\textit{Formula 27.}
$$
[x_{-a-b}(u),x_{2a+b}(t)]=
$$
$$
=x_a(C_{11,2a+b,-a-b}(-t)u)
 x_{3a+b}(C_{21,2a+b,-a-b}(-t)^2u)\,x_{-b}(C_{12,2a+b,-a-b}(-t)u^2).
$$
Calculation of constants: \vskip3pt

$\displaystyle{C_{11,2a+b,-a-b}=M_{2a+b,-a-b,1}=N_{2a+b,-a-b}=2\epsilon_2},$
\vskip7pt

$\displaystyle{C_{21,2a+b,-a-b}=M_{2a+b,-a-b,2}=\frac{1}{2!}N_{2a+b,-a-b}N_{2a+b,a}=\frac{1}{2}(2\epsilon_2)(-3\epsilon_3)=-3\epsilon_2\epsilon_3
,}$ \vskip7pt

$\displaystyle{C_{12,2a+b,-a-b}=(-1)^2M_{-a-b,2a+b,2}=\frac{1}{2!}N_{-a-b,2a+b}N_{-a-b,a}=\frac{1}{2}(-2\epsilon_2)(3\epsilon_1)=
-3\epsilon_1\epsilon_2.}$ \vskip7pt

\noindent \textbf{Final formula:}
\vskip10pt

\begin{center}
 \fbox{$\displaystyle{[x_{-a-b}(u),x_{2a+b}(t)]=x_a(-2\epsilon_2tu)\,x_{3a+b}(-3\epsilon_2\epsilon_3t^2u)\,x_{-b}(3\epsilon_1\epsilon_2tu^2).}$}
\end{center}
\bigskip

\noindent\textit{Formula 28.}
$$
[x_{-3a-b}(u),x_{2a+b}(t)]=x_{-a}(C_{11,2a+b,-3a-b}(-t)u)\times
$$
$$
x_{a+b}(C_{21,2a+b,-3a-b}(-t)^2u)\,x_{3a+2b}(C_{31,2a+b,-3a-b}(-t)^3u)\,x_b(C_{32,2a+b,-3a-b}(-t)^3u^2).
$$
Calculation of constants: \vskip3pt

$\displaystyle{C_{11,2a+b,-3a-b}=M_{2a+b,-3a-b,1}=N_{2a+b,-3a-b}=\epsilon_3},$
\vskip7pt

$\displaystyle{C_{21,2a+b,-3a-b}=M_{2a+b,-3a-b,2}=\frac{1}{2!}N_{2a+b,-3a-b}N_{2a+b,-a}=\frac{1}{2}(\epsilon_3)(-2\epsilon_2)=-\epsilon_2\epsilon_3
,}$ \vskip7pt

$\displaystyle{C_{31,2a+b,-3a-b}=M_{2a+b,-3a-b,3}=\frac{1}{3!}N_{2a+b,-3a-b}N_{2a+b,-a}N_{2a+b,a+b}=\frac{1}{6}(\epsilon_3)(-2\epsilon_2)(3\epsilon_5)=}$

$\displaystyle{\ \ \ \ \ \ \ \ \ \ \ \ \ \ \ \ \ \
=-\epsilon_2\epsilon_3\epsilon_5,}$ \vskip7pt

$\displaystyle{C_{32,2a+b,-3a-b}=\frac{1}{3}M_{-a,2a+b,2}=\frac{1}{3!}N_{-a,2a+b}N_{-a,a+b}=\frac{1}{6}(2\epsilon_2)(3\epsilon_1)=
\epsilon_1\epsilon_2.}$ \vskip7pt

\noindent \textbf{Final formula:}
\vskip10pt

\begin{center}
 \fbox{$\displaystyle{[x_{-3a-b}(u),x_{2a+b}(t)]=x_{-a}(-\epsilon_3tu)\,x_{a+b}(-\epsilon_2\epsilon_3t^2u)\,x_{3a+2b}(\epsilon_2\epsilon_3\epsilon_5t^3u)\,
x_b(-\epsilon_1\epsilon_2t^3u^2).}$}
\end{center}
\bigskip

\noindent\textit{Formula 29.}
$$
[x_{2a+b}(u),x_{-3a-b}(t)]=x_{-a}(C_{11,-3a-b,2a+b}(-t)u)\times
$$
$$
x_{a+b}(C_{12,-3a-b,2a+b}(-t)u^2)\,x_{3a+2b}(C_{13,-3a-b,2a+b}(-t)u^2)\,x_b(C_{23,-3a-b,2a+b}(-t)^2u^3).
$$
Calculation of constants: \vskip3pt

$\displaystyle{C_{11,-3a-b,2a+b}=M_{-3a-b,2a+b,1}=N_{-3a-b,2a+b}=-\epsilon_3},$
\vskip7pt

$\displaystyle{C_{12,-3a-b,2a+b}=(-1)^2M_{2a+b,-3a-b,2}=\frac{1}{2!}N_{2a+b,-3a-b}N_{2a+b,-a}=\frac{1}{2}(\epsilon_3)(-2\epsilon_2)=-\epsilon_2\epsilon_3
,}$ \vskip7pt

$\displaystyle{C_{13,-3a-b,2a+b}=(-1)^3M_{2a+b,-3a-b,3}=-\frac{1}{3!}N_{2a+b,-3a-b}N_{2a+b,-a}N_{2a+b,a+b}}$\vskip7pt

$\displaystyle{\ \ \ \ \ \ \ \ \ \ \ \ \ \ \ \ \ \
=-\frac{1}{6}(\epsilon_3)(-2\epsilon_2)(3\epsilon_5)=
\epsilon_2\epsilon_3\epsilon_5,}$ \vskip7pt

$\displaystyle{C_{23,-3a-b,2a+b}=-\frac{2}{3}M_{-a,2a+b,2}=-\frac{1}{3}N_{-a,2a+b}N_{-a,a+b}=-\frac{1}{3}(2\epsilon_2)(3\epsilon_1)=
-2\epsilon_1\epsilon_2.}$ \vskip7pt

\noindent \textbf{Final formula:}
\vskip10pt

\begin{center}
 \fbox{$\displaystyle{[x_{2a+b}(u),x_{-3a-b}(t)]=x_{-a}(\epsilon_3tu)\,x_{a+b}(\epsilon_2\epsilon_3tu^2)\,x_{3a+2b}(-\epsilon_2\epsilon_3\epsilon_5tu^3)\,
x_b(-2\epsilon_1\epsilon_2t^2u^3).}$}
\end{center}
\bigskip

\noindent\textit{Formula 30.}
$$
[x_{-3a-2b}(u),x_{2a+b}(t)]=x_{-a-b}(C_{11,2a+b,-3a-2b}(-t)u)\times
$$
$$
x_{a}(C_{21,2a+b,-3a-2b}(-t)^2u)\,x_{3a+b}(C_{31,2a+b,-3a-2b}(-t)^3u)\,x_{-b}(C_{32,2a+b,-3a-2b}(-t)^3u^2).
$$
Calculation of constants: \vskip3pt

$\displaystyle{C_{11,2a+b,-3a-2b}=M_{2a+b,-3a-2b,1}=N_{2a+b,-3a-2b}=-\epsilon_5},$
\vskip7pt

$\displaystyle{C_{21,2a+b,-3a-2b}=M_{2a+b,-3a-2b,2}=\frac{1}{2!}N_{2a+b,-3a-2b}N_{2a+b,-a-b}=\frac{1}{2}(-\epsilon_5)(2\epsilon_2)=-\epsilon_2\epsilon_5
,}$ \vskip7pt

$\displaystyle{C_{31,2a+b,-3a-2b}=M_{2a+b,-3a-2b,3}=\frac{1}{3!}N_{2a+b,-3a-2b}N_{2a+b,-a-b}N_{2a+b,a}=}$
\vskip7pt

$\qquad\qquad\qquad\
\displaystyle{=\frac{1}{6}(-\epsilon_5)(2\epsilon_2)(-3\epsilon_3)=\epsilon_2\epsilon_3\epsilon_5,}$ \vskip7pt

$\displaystyle{C_{32,2a+b,-3a-2b}=\frac{1}{3}M_{-a-b,2a+b,2}=\frac{1}{3!}N_{-a-b,2a+b}N_{-a-b,a}=\frac{1}{6}(-2\epsilon_2)(3\epsilon_1)=
-\epsilon_1\epsilon_2.}$ \vskip7pt

\noindent \textbf{Final formula:}
\vskip10pt

\begin{center}
 \fbox{$\displaystyle{[x_{-3a-2b}(u),x_{2a+b}(t)]=x_{-a-b}(\epsilon_5tu)\,x_a(-\epsilon_2\epsilon_5t^2u)\,x_{3a+b}(-\epsilon_2\epsilon_3\epsilon_5t^3u)\,
x_{-b}(\epsilon_1\epsilon_2t^3u^2).}$}
\end{center}
\bigskip

\noindent\textit{Formula 31.}
$$
[x_{2a+b}(u),x_{-3a-2b}(t)]=x_{-a-b}(C_{11,-3a-2b,2a+b}(-t)u)\times
$$
$$
x_a(C_{12,-3a-2b,2a+b}(-t)u^2)\,x_{3a+b}(C_{13,-3a-2b,2a+b}(-t)u^3)\,x_{-b}(C_{23,-3a-2b,2a+b}(-t)^2u^3).
$$
Calculation of constants: \vskip3pt

$\displaystyle{C_{11,-3a-2b,2a+b}=M_{-3a-2b,2a+b,1}=N_{-3a-2b,2a+b}=\epsilon_5},$
\vskip7pt

$\displaystyle{C_{12,-3a-2b,2a+b}=(-1)^2M_{2a+b,-3a-2b,2}=\frac{1}{2!}N_{2a+b,-3a-2b}N_{2a+b,-a-b}=\frac{1}{2}(-\epsilon_5)(2\epsilon_2)=}$
\vskip7pt

$\displaystyle{\ \ \ \ \ \ \ \ \ \ \ \ \ \ \ \ \ \ \
=-\epsilon_2\epsilon_5,}$ \vskip7pt

$\displaystyle{C_{13,-3a-2b,2a+b}=(-1)^3M_{2a+b,-3a-2b,3}=-\frac{1}{3!}N_{2a+b,-3a-2b}N_{2a+b,-a-b}N_{2a+b,a}=}$\vskip7pt

$\displaystyle{\ \ \ \ \ \ \ \ \ \ \ \ \ \ \ \ \ \ \
=-\frac{1}{6}(-\epsilon_5)(2\epsilon_2)(-3\epsilon_3)=-\epsilon_2\epsilon_3\epsilon_5,}$
\vskip7pt

$\displaystyle{C_{23,-3a-2b,2a+b}=-\frac{2}{3}M_{-a-b,2a+b,2}=-\frac{1}{3}N_{-a-b,2a+b}N_{-a-b,a}=-\frac{1}{3}(-2\epsilon_2)(3\epsilon_1)=
2\epsilon_1\epsilon_2.}$ \vskip7pt

\noindent \textbf{Final formula:}
\vskip10pt

\begin{center}
 \fbox{$\displaystyle{[x_{2a+b}(u),x_{-3a-2b}(t)]=x_{-a-b}(-\epsilon_5tu)\,x_a(\epsilon_2\epsilon_5tu^2)\,x_{3a+b}(\epsilon_2\epsilon_3\epsilon_5tu^3)\,
x_{-b}(2\epsilon_1\epsilon_2t^2u^3).}$}
\end{center}
\bigskip

\noindent\textit{Formula 32.}
$$
[x_{-a}(u),x_{3a+b}(t)]=x_{2a+b}(C_{11,3a+b,-a}(-t)u)\times
$$
$$
x_{a+b}(C_{12,3a+b,-a}(-t)u^2)\, x_{b}(C_{13,3a+b,-a}(-t)u^3)\,x_{3a+2b}(C_{23,3a+b,-a}(-t)^2u^3).
$$
Calculation of constants: \vskip3pt

 $\displaystyle{C_{11,3a+b,-a}=M_{3a+b,-a,1}=N_{3a+b,-a}=-\epsilon_3},$
\vskip7pt

 $\displaystyle{C_{12,3a+b,-a}=(-1)^2M_{-a,3a+b,2}=\frac{1}{2!}N_{-a,3a+b}N_{-a,2a+b}=
 \frac{1}{2}(\epsilon_3)(2\epsilon_2)=\epsilon_2\epsilon_3
,}$ \vskip7pt

 $\displaystyle{C_{13,3a+b,-a}=(-1)^3M_{-a,3a+b,3}=-\frac{1}{3!}N_{-a,3a+b}N_{-a,2a+b}N_{-a,a+b}=
 -\frac{1}{6}(\epsilon_3)(2\epsilon_2)(3\epsilon_1)=}$
\vskip7pt

$\displaystyle{\ \ \ \ \ \ \ \ \ \ \ \ \ \ \,=-\epsilon_1\epsilon_2\epsilon_3,}$ 
\vskip7pt

 $\displaystyle{C_{23,3a+b,-a}=-\frac23 M_{2a+b,-a,2}=-\frac{1}{3}N_{2a+b,-a}N_{2a+b,a+b}=
 -\frac{1}{3}(-2\epsilon_2)(3\epsilon_5)=2\epsilon_2\epsilon_5
.}$ \vskip7pt

\noindent \textbf{Final formula:}
\vskip10pt

\begin{center}
 \fbox{$\displaystyle{[x_{-a}(u),x_{3a+b}(t)]=
 x_{2a+b}(\epsilon_3tu)\,x_{a+b}(-\epsilon_2\epsilon_3 tu^2)\,x_b(\epsilon_1\epsilon_2\epsilon_3 tu^3)\,
 x_{3a+2b}(2\epsilon_2\epsilon_5t^2u^3).}$}
\end{center}
\bigskip

\noindent\textit{Formula 33.}
$$
[x_{3a+b}(u),x_{-a}(t)]=x_{2a+b}(C_{11,-a,3a+b}(-t)u)\times
$$
$$
x_{a+b}(C_{21,-a,3a+b}(-t)^2u)\, x_b(C_{31,-a,3a+b}(-t)^3u)\,x_{3a+2b}(C_{32,-a,3a+b}(-t)^3u^2).
$$
Calculation of constants: \vskip3pt

$\displaystyle{C_{11,-a,3a+b}=M_{-a,3a+b,1}=N_{-a,3a+b}=\epsilon_3},$
\vskip7pt

$\displaystyle{C_{21,-a,3a+b}=M_{-a,3a+b,2}=\frac{1}{2!}N_{-a,3a+b}N_{-a,2a+b}=
 \frac{1}{2}(\epsilon_3)(2\epsilon_2)=\epsilon_2\epsilon_3
,}$ \vskip7pt

$\displaystyle{C_{31,-a,3a+b}=M_{-a,3a+b,3}=\frac{1}{3!}N_{-a,3a+b}N_{-a,2a+b}N_{-a,a+b}=
 \frac{1}{6}(\epsilon_3)(2\epsilon_2)(3\epsilon_1)=\epsilon_1\epsilon_2\epsilon_3
,}$ \vskip7pt

$\displaystyle{C_{32,-a,3a+b}=\frac{1}{3}M_{2a+b,-a,2}=\frac{1}{3!}N_{2a+b,-a}N_{2a+b,a+b}=
 \frac{1}{6}(-2\epsilon_2)(3\epsilon_5)=-\epsilon_2\epsilon_5
.}$ \vskip7pt

\noindent \textbf{Final formula:}
\vskip10pt

\begin{center}
 \fbox{$\displaystyle{[x_{3a+b}(u),x_{-a}(t)]=
 x_{2a+b}(-\epsilon_3tu)\,x_{a+b}(\epsilon_2\epsilon_3 t^2u)\,x_b(-\epsilon_1\epsilon_2\epsilon_3 t^3u)\,x_{3a+2b}(\epsilon_2\epsilon_5t^3u^2).}$}
\end{center}
\bigskip

\noindent\textit{Formula 34.}
$$
[x_{-2a-b}(u),x_{3a+b}(t)]=x_a(C_{11,3a+b,-2a-b}(-t)u)\times
$$
$$
 x_{-a-b}(C_{12,3a+b,-2a-b}(-t)u^2)\, x_{-3a-2b}(C_{13,3a+b,-2a-b}(-t)u^3)\,x_{-b}(C_{23,3a+b,-2a-b}(-t)^2u^3).
$$
Calculation of constants: \vskip3pt

 $\displaystyle{C_{11,3a+b,-2a-b}=M_{3a+b,-2a-b,1}=N_{3a+b,-2a-b}=\epsilon_3},$
\vskip7pt

 $\displaystyle{C_{12,3a+b,-2a-b}=(-1)^2M_{-2a-b,3a+b,2}=\frac{1}{2!}N_{-2a-b,3a+b}N_{-2a-b,a}=
 \frac{1}{2}(-\epsilon_3)(2\epsilon_2)=-\epsilon_2\epsilon_3
,}$ \vskip7pt

 $\displaystyle{C_{13,3a+b,-2a-b}=(-1)^3M_{-2a-b,3a+b,3}=-\frac{1}{3!}N_{-2a-b,3a+b}N_{-2a-b,a}N_{-2a-b,-a-b}=}$

 $\displaystyle{\ \ \ \ \ \ \ \ \ \ \ \ \ \ \ \ \ \ =
 -\frac{1}{6}(-\epsilon_3)(2\epsilon_2)(-3\epsilon_5)=-\epsilon_2\epsilon_3\epsilon_5,}$ 
\vskip7pt

 $\displaystyle{C_{23,3a+b,-2a-b}=-\frac23 M_{a,-2a-b,2}=-\frac{1}{3}N_{a,-2a-b}N_{a,-a-b}=
 -\frac{1}{3}(-2\epsilon_2)(-3\epsilon_1)=-2\epsilon_1\epsilon_2
.}$ \vskip7pt

\noindent \textbf{Final formula:}
\vskip10pt

\begin{center}
 \fbox{$\displaystyle{[x_{-2a-b}(u),x_{3a+b}(t)]=
 x_a(-\epsilon_3tu)\,x_{-a-b}(\epsilon_2\epsilon_3 tu^2)\,x_{-3a-2b}(\epsilon_2\epsilon_3\epsilon_5 tu^3)\,
 x_{-b}(-2\epsilon_1\epsilon_2t^2u^3).}$}
\end{center}
\bigskip

\noindent\textit{Formula 35.}
$$
[x_{3a+b}(u),x_{-2a-b}(t)]=x_a(C_{11,-2a-b,3a+b}(-t)u)\times
$$
$$
 x_{-a-b}(C_{21,-2a-b,3a+b}(-t)^2u)\, x_{-3a-2b}(C_{31,-2a-b,3a+b}(-t)^3u)\,x_{-b}(C_{32,-2a-b,3a+b}(-t)^3u^2).
$$
Calculation of constants: \vskip3pt

$\displaystyle{C_{11,-2a-b,3a+b}=M_{-2a-b,3a+b,1}=N_{-2a-b,3a+b}=-\epsilon_3},$
\vskip7pt

$\displaystyle{C_{21,-2a-b,3a+b}=M_{-2a-b,3a+b,2}=\frac{1}{2!}N_{-2a-b,3a+b}N_{-2a-b,a}=
 \frac{1}{2}(-\epsilon_3)(2\epsilon_2)=-\epsilon_2\epsilon_3
,}$ \vskip7pt

$\displaystyle{C_{31,-2a-b,3a+b}=M_{-2a-b,3a+b,3}=\frac{1}{3!}N_{-2a-b,3a+b}N_{-2a-b,a}N_{-2a-b,-a-b}=}$

$\qquad\qquad\qquad \displaystyle{=\frac{1}{6}(-\epsilon_3)(2\epsilon_2)(-3\epsilon_5)=\epsilon_2\epsilon_3\epsilon_5
,}$ \vskip7pt

$\displaystyle{C_{32,-2a-b,3a+b}=\frac{1}{3}M_{a,-2a-b,2}=\frac{1}{3!}N_{a,-2a-b}N_{a,-a-b}=
 \frac{1}{6}(-2\epsilon_2)(-3\epsilon_1)=\epsilon_1\epsilon_2
.}$ \vskip7pt

\noindent \textbf{Final formula:}
\vskip10pt

\begin{center}
 \fbox{$\displaystyle{[x_{3a+b}(u),x_{-2a-b}(t)]=
 x_a(\epsilon_3tu)\,x_{-a-b}(-\epsilon_2\epsilon_3 t^2u)\,x_{-3a-2b}(-\epsilon_2\epsilon_3\epsilon_5 t^3u)\,x_{-b}(-\epsilon_1\epsilon_2t^3u^2).}$}
\end{center}
\bigskip

\noindent\textit{Formula 36.}
$$
[x_{-3a-2b}(u),x_{3a+b}(t)]=x_{-b}(C_{11,3a+b,-3a-2b}(-t)u).
$$

Calculation of constants: \vskip3pt

$\displaystyle{C_{11,3a+b,-3a-2b}=M_{3a+b,-3a-2b,1}=N_{3a+b,-3a-2b}=\epsilon_4}.$
\vskip7pt

\noindent \textbf{Final formula:}
\vskip10pt

\begin{center}
 \fbox{$\displaystyle{[x_{-3a-2b}(u),x_{3a+b}(t)]=x_{-b}(-\epsilon_4tu).}$}
\end{center}
\bigskip

\noindent\textit{Formula 37.}
$$
[x_{-a-b}(u),x_{3a+2b}(t)]=x_{2a+b}(C_{11,3a+2b,-a-b}(-t)u)\times
$$
$$
 x_a(C_{12,3a+2b,-a-b}(-t)u^2)\, x_{-b}(C_{13,3a+2b,-a-b}(-t)u^3)\,x_{3a+b}(C_{23,3a+2b,-a-b}(-t)^2u^3).
$$
Calculation of constants: \vskip3pt

 $\displaystyle{C_{11,3a+2b,-a-b}=M_{3a+2b,-a-b,1}=N_{3a+2b,-a-b}=\epsilon_5},$
\vskip7pt

 $\displaystyle{C_{12,3a+2b,-a-b}=(-1)^2M_{-a-b,3a+2b,2}=\frac{1}{2!}N_{-a-b,3a+2b}N_{-a-b,2a+b}=
 \frac{1}{2}(-\epsilon_5)(-2\epsilon_2)=\epsilon_2\epsilon_5
,}$ \vskip7pt

 $\displaystyle{C_{13,3a+2b,-a-b}=(-1)^3M_{-a-b,3a+2b,3}=-\frac{1}{3!}N_{-a-b,3a+2b}N_{-a-b,2a+b}N_{-a-b,a}=}$

 $\displaystyle{\ \ \ \ \ \ \ \ \ \ \ \ \ \ \ \ \ \ =
 -\frac{1}{6}(-\epsilon_5)(-2\epsilon_2)(3\epsilon_1)=-\epsilon_1\epsilon_2\epsilon_5
,}$ \vskip7pt

 $\displaystyle{C_{23,3a+2b,-a-b}=-\frac23 M_{2a+b,-a-b,2}=-\frac{1}{3}N_{2a+b,-a-b}N_{2a+b,a}=
 -\frac{1}{3}(2\epsilon_2)(-3\epsilon_3)=2\epsilon_2\epsilon_3
.}$ \vskip7pt

\noindent \textbf{Final formula:}
\vskip10pt

\begin{center}
 \fbox{$\displaystyle{[x_{-a-b}(u),x_{3a+2b}(t)]=
 x_{2a+b}(-\epsilon_5tu)\,x_a(-\epsilon_2\epsilon_5 tu^2)\,x_{-b}(\epsilon_1\epsilon_2\epsilon_5 tu^3)\,
 x_{3a+b}(2\epsilon_2\epsilon_3t^2u^3).}$}
\end{center}
\bigskip

\noindent\textit{Formula 38.}
$$
[x_{3a+2b}(u),x_{-a-b}(t)]=x_{2a+b}(C_{11,-a-b,3a+2b}(-t)u)\times
$$
$$
 x_a(C_{21,-a-b,3a+2b}(-t)^2u)\, x_{-b}(C_{31,-a-b,3a+2b}(-t)^3u)\,x_{3a+b}(C_{32,-a-b,3a+2b}(-t)^3u^2).
$$
Calculation of constants: \vskip3pt

$\displaystyle{C_{11,-a-b,3a+2b}=M_{-a-b,3a+2b,1}=N_{-a-b,3a+2b}=-\epsilon_5},$
\vskip7pt

$\displaystyle{C_{21,-a-b,3a+2b}=M_{-a-b,3a+2b,2}=\frac{1}{2!}N_{-a-b,3a+2b}N_{-a-b,2a+b}=
 \frac{1}{2}(-\epsilon_5)(-2\epsilon_2)=\epsilon_2\epsilon_5
,}$ \vskip7pt

$\displaystyle{C_{31,-a-b,3a+2b}=M_{-a-b,3a+2b,3}=\frac{1}{3!}N_{-a-b,3a+2b}N_{-a-b,2a+b}N_{-a-b,a}=}$

$\qquad\qquad\qquad \displaystyle{=\frac{1}{6}(-\epsilon_5)(-2\epsilon_2)(3\epsilon_1)=\epsilon_1\epsilon_2\epsilon_5,}$ \vskip7pt

$\displaystyle{C_{32,-a-b,3a+2b}=\frac{1}{3}M_{2a+b,-a-b,2}=\frac{1}{3!}N_{2a+b,-a-b}N_{2a+b,a}=
 \frac{1}{6}(2\epsilon_2)(-3\epsilon_3)=-\epsilon_2\epsilon_3
.}$ \vskip7pt

\noindent \textbf{Final formula:}
\vskip10pt

\begin{center}
 \fbox{$\displaystyle{[x_{3a+2b}(u),x_{-a-b}(t)]=
 x_{2a+b}(\epsilon_5tu)\,x_a(\epsilon_2\epsilon_5 t^2u)\,x_{-b}(-\epsilon_1\epsilon_2\epsilon_5 t^3u)\,x_{3a+b}(\epsilon_2\epsilon_3t^3u^2).}$}
\end{center}
\bigskip

\noindent\textit{Formula 39.}
$$
[x_{-2a-b}(u),x_{3a+2b}(t)]=x_{a+b}(C_{11,3a+2b,-2a-b}(-t)u)\times
$$
$$
 x_{-a}(C_{12,3a+2b,-2a-b}(-t)u^2)\, x_{-3a-b}(C_{13,3a+2b,-2a-b}(-t)u^3)\,x_b(C_{23,3a+2b,-2a-b}(-t)^2u^3).
$$
Calculation of constants: \vskip3pt

 $\displaystyle{C_{11,3a+2b,-2a-b}=M_{3a+2b,-2a-b,1}=N_{3a+2b,-2a-b}=-\epsilon_5},$
\vskip7pt

 $\displaystyle{C_{12,3a+2b,-2a-b}=(-1)^2M_{-2a-b,3a+2b,2}=\frac{1}{2!}N_{-2a-b,3a+2b}N_{-2a-b,a+b}=
 \frac{1}{2}(\epsilon_5)(-2\epsilon_2)=}$
\vskip7pt

$\displaystyle{\ \ \ \ \ \ \ \ \ \ \ \ \ \ \ \ \ \ \ = -\epsilon_2\epsilon_5,}$ \vskip7pt

 $\displaystyle{C_{13,3a+2b,-2a-b}=(-1)^3M_{-2a-b,3a+2b,3}=-\frac{1}{3!}N_{-2a-b,3a+2b}N_{-2a-b,a+b}N_{-2a-b,-a}=}$

 $\displaystyle{\ \ \ \ \ \ \ \ \ \ \ \ \ \ \ \ \ \ \ =
 -\frac{1}{6}(\epsilon_5)(-2\epsilon_2)(3\epsilon_3)=\epsilon_2\epsilon_3\epsilon_5
,}$ \vskip7pt

 $\displaystyle{C_{23,3a+2b,-2a-b}=-\frac23 M_{a+b,-2a-b,2}=-\frac{1}{3}N_{a+b,-2a-b}N_{a+b,-a}=
 -\frac{1}{3}(2\epsilon_2)(-3\epsilon_1)=2\epsilon_1\epsilon_2
.}$ \vskip7pt

\noindent \textbf{Final formula:}
\vskip10pt

\begin{center}
 \fbox{$\displaystyle{[x_{-2a-b}(u),x_{3a+2b}(t)]=
 x_{a+b}(\epsilon_5tu)\,x_{-a}(\epsilon_2\epsilon_5 tu^2)\,x_{-3a-b}(-\epsilon_2\epsilon_3\epsilon_5 tu^3)\,
 x_b(2\epsilon_1\epsilon_2t^2u^3).}$}
\end{center}
\bigskip

\noindent\textit{Formula 40.}
$$
[x_{3a+2b}(u),x_{-2a-b}(t)]=x_{a+b}(C_{11,-2a-b,3a+2b}(-t)u)\times
$$
$$
 x_{-a}(C_{21,-2a-b,3a+2b}(-t)^2u)\, x_{-3a-b}(C_{31,-2a-b,3a+2b}(-t)^3u)\,x_b(C_{32,-2a-b,3a+2b}(-t)^3u^2).
$$
Calculation of constants: \vskip3pt

$\displaystyle{C_{11,-2a-b,3a+2b}=M_{-2a-b,3a+2b,1}=N_{-2a-b,3a+2b}=\epsilon_5},$
\vskip7pt

$\displaystyle{C_{21,-2a-b,3a+2b}=M_{-2a-b,3a+2b,2}=\frac{1}{2!}N_{-2a-b,3a+2b}N_{-2a-b,a+b}=
 \frac{1}{2}(\epsilon_5)(-2\epsilon_2)=-\epsilon_2\epsilon_5
,}$ \vskip7pt

$\displaystyle{C_{31,-2a-b,3a+2b}=M_{-2a-b,3a+2b,3}=\frac{1}{3!}N_{-2a-b,3a+2b}N_{-2a-b,a+b}N_{-2a-b,-a}=}$
\vskip7pt

$\qquad\qquad\qquad\ =\displaystyle{=\frac{1}{6}(\epsilon_5)(-2\epsilon_2)(3\epsilon_3)=-\epsilon_2\epsilon_3\epsilon_5 ,}$ \vskip7pt

$\displaystyle{C_{32,-2a-b,3a+2b}=\frac{1}{3}M_{a+b,-2a-b,2}=\frac{1}{3!}N_{a+b,-2a-b}N_{a+b,-a}=
 \frac{1}{6}(2\epsilon_2)(-3\epsilon_1)=-\epsilon_1\epsilon_2
.}$ \vskip7pt

\noindent \textbf{Final formula:}
\vskip10pt

\begin{center}
 \fbox{$\displaystyle{[x_{3a+2b}(u),x_{-2a-b}(t)]=
 x_{a+b}(-\epsilon_5tu)\,x_{-a}(-\epsilon_2\epsilon_5 t^2u)\,x_{-3a-b}(\epsilon_2\epsilon_3\epsilon_5 t^3u)\,x_b(\epsilon_1\epsilon_2t^3u^2).}$}
\end{center}
\bigskip

\noindent\textit{Formula 41.}
$$
[x_{-3a-b}(u),x_{3a+2b}(t)]=x_b(C_{11,3a+2b,-3a-b}(-t)u).
$$

Calculation of constants: \vskip3pt

$\displaystyle{C_{11,3a+2b,-3a-b}=M_{3a+2b,-3a-b,1}=N_{3a+2b,-3a-b}=\epsilon_4}.$
\vskip7pt

\noindent \textbf{Final formula:}
\vskip10pt

\begin{center}
 \fbox{$\displaystyle{[x_{-3a-b}(u),x_{3a+2b}(t)]=x_{a+b}(-\epsilon_4tu).}$}
\end{center}

\section{Special case}

\subsection{Table of the Structure Constants}

Using Theorem 1 we obtain

\textbf{Theorem 2.} \textit{Let the values of all structure constants corresponding to extraspecial pairs be positive.
Then all structure constants have the values indicated in the following table.}
\vskip10pt

{\scriptsize
\begin{center}
\begin{tabular}{|r|r|r|r|r|r|r|r|r|r|r|r|r|}
\hline
 $N_{r,s}$ & $-3a-2b$ & $-3a-b$ & $-2a-b$ & $-a-b$ & $-a$ & $-b$ & $a$ & $b$ & $a+b$ & $2a+b$ & $3a+b$ & $3a+2b$ \\ \hline
 $-3a-2b$  &          &         &         &        &      &      &     &   1 &    -1 &      1 &     -1 &         \\ \hline
 $-3a-b$   &          &         &         &        &      &    1 &   1 &     &       &     -1 &        &      -1 \\ \hline
 $-2a-b$   &          &         &         &     -3 &    3 &      &   2 &     &    -2 &        &     -1 &       1 \\ \hline
 $-a-b$    &          &         &       3 &        &    2 &      &   3 &  -1 &       &     -2 &        &      -1 \\ \hline
 $-b$      &          &      -1 &         &        &    1 &      &     &     &    -1 &        &        &       1 \\ \hline
 $-a$      &          &         &      -3 &     -2 &      &   -1 &     &     &     3 &      2 &      1 &         \\ \hline
 $a$       &          &      -1 &      -2 &     -3 &      &      &     &   1 &     2 &      3 &        &         \\ \hline
 $b$       &       -1 &         &         &      1 &      &      &  -1 &     &       &        &      1 &         \\ \hline
 $a+b$     &        1 &         &       2 &        &   -3 &    1 &  -2 &     &       &     -3 &        &         \\ \hline
 $2a+b$    &       -1 &       1 &         &      2 &   -2 &      &  -3 &     &     3 &        &        &         \\ \hline
 $3a+b$    &        1 &         &       1 &        &   -1 &      &     &  -1 &       &        &        &         \\ \hline
 $3a+2b$   &          &       1 &      -1 &      1 &      &   -1 &     &     &       &        &        &         \\ \hline
\end{tabular}
\end{center}
}
\medskip

\centerline{\textbf{Table 2.}}

\vskip5mm

\subsection{List of Formulas}

\centerline{\textbf{Positive roots}} \vskip3pt

\begin{equation}\label{pf1}
[x_b(u),x_a(t)]
 =x_{a+b}(-tu)\,x_{2a+b}(t^2u)\,x_{3a+b}(-t^3u)\,x_{3a+2b}(-t^3u^2),
\end{equation}

\begin{equation}\label{pf2}
[x_a(u),x_b(t)]
 =x_{a+b}(tu)\,x_{2a+b}(-tu^2)\,x_{3a+b}(tu^3)\,x_{3a+2b}(-2t^2u^3),
\end{equation}

\begin{equation}\label{pf3}
[x_{a+b}(u),x_a(t)]
 =x_{2a+b}(-2tu)\,x_{3a+b}(3t^2u)\,x_{3a+2b}(-3tu^2),
\end{equation}

\begin{equation}\label{pf4}
[x_{2a+b}(u),x_a(t)]=x_{3a+b}(-3tu),
\end{equation}

\begin{equation}\label{pf5}
[x_{3a+b}(u),x_b(t)]=x_{3a+2b}(-tu),
\end{equation}

\begin{equation}\label{pf6}
[x_{2a+b}(u),x_{a+b}(t)]=x_{3a+2b}(3tu).
\end{equation}

\vskip10pt \centerline{\textbf{Negative roots}} \vskip5pt

\begin{equation}\label{pf7}
[x_{-b}(u),x_{-a}(t)]
 =x_{-a-b}(tu)\,x_{-2a-b}(t^2u)\,x_{-3a-b}(t^3u)\,x_{-3a-2b}(-t^3u^2),
\end{equation}

\begin{equation}\label{pf8}
[x_{-a}(u),x_{-b}(t)]
 =x_{-a-b}(-tu)\,x_{-2a-b}(-tu^2)\,x_{-3a-b}(-tu^3)\,x_{-3a-2b}(-2t^2u^3),
\end{equation}

\begin{equation}\label{pf9}
[x_{-a-b}(u),x_{-a}(t)]=x_{-2a-b}(2tu)\,x_{-3a-b}(3t^2u)\,x_{-3a-2b}(-3tu^2),
\end{equation}

\begin{equation}\label{pf10}
[x_{-2a-b}(u),x_{-a}(t)]=x_{-3a-b}(3tu),
\end{equation}

\begin{equation}\label{pf11}
[x_{-3a-b}(u),x_{-b}(t)]=x_{-3a-2b}(tu),
\end{equation}

\begin{equation}\label{pf12}
[x_{-2a-b}(u),x_{-a-b}(t)]=x_{-3a-2b}(-3tu).
\end{equation}

\vskip10pt \centerline{\textbf{Positive and negative roots}} \vskip5pt

\begin{equation}\label{pf13}
[x_{-a-b}(u),x_{a}(t)]=x_{-b}(3tu),
\end{equation}

\begin{equation}\label{pf14}
[x_{-2a-b}(u),x_a(t)]=x_{-a-b}(2tu)\,x_{-b}(3t^2u)\,x_{-3a-2b}(3tu^2),
\end{equation}

\begin{equation}\label{pf15}
[x_{-3a-b}(u),x_a(t)]
 =x_{-2a-b}(tu)\,x_{-a-b}(t^2u)\,x_{-b}(t^3u)\,x_{-3a-2b}(t^3u^2),
\end{equation}

\begin{equation}\label{pf16}
[x_{a}(u),x_{-3a-b}(t)]
 =x_{-2a-b}(-tu)\,x_{-a-b}(-tu^2)\,x_{-b}(-tu^3)\,x_{-3a-2b}(2t^2u^3),
\end{equation}

\begin{equation}\label{pf17}
[x_b(u),x_{-a-b}(t)]
 =x_{-a}(tu)\,x_{-2a-b}(-t^2u)\,x_{-3a-2b}(t^3u)\,x_{-3a-b}(-t^3u^2),
\end{equation}

\begin{equation}\label{pf18}
[x_{-a-b}(u),x_b(t)]
 =x_{-a}(-tu)\,x_{-2a-b}(tu^2)\,x_{-3a-2b}(-tu^3)\,x_{-3a-b}(-2t^2u^3),
\end{equation}

\begin{equation}\label{pf19}
[x_{-3a-2b}(u),x_b(t)]=x_{-3a-b}(tu),
\end{equation}

\begin{equation}\label{pf20}
[x_{-a}(u),x_{a+b}(t)]=x_b(3tu),
\end{equation}

\begin{equation}\label{pf21}
[x_{-b}(u),x_{a+b}(t)]
 =x_a(-tu)\,x_{2a+b}(-t^2u)\,x_{3a+2b}(-t^3u)\,x_{3a+b}(-t^3u^2),
\end{equation}

\begin{equation}\label{pf22}
[x_{a+b}(u),x_{-b}(t)]
 =x_a(tu)\,x_{2a+b}(tu^2)\,x_{3a+2b}(tu^3)\, x_{3a+b}(-2t^2u^3),
\end{equation}

\begin{equation}\label{pf23}
[x_{-2a-b}(u),x_{a+b}(t)]=
 x_{-a}(-2tu)\,x_b(-3t^2u)\,x_{-3a-b}(3tu^2),
\end{equation}

\begin{equation}\label{pf24}
[x_{-3a-2b}(u),x_{a+b}(t)]
 =x_{-2a-b}(-tu)\,x_{-a}(t^2u)\,x_b(t^3u)\, x_{-3a-b}(t^3u^2),
\end{equation}

\begin{equation}\label{pf25}
[x_{a+b}(u),x_{-3a-2b}(t)]
 =x_{-2a-b}(tu)\,x_{-a}(-tu^2)\,x_b(-tu^3)\, x_{-3a-b}(2t^2u^3),
\end{equation}

\begin{equation}\label{pf26}
[x_{-a}(u),x_{2a+b}(t)]=
 x_{a+b}(2tu)\,x_{3a+2b}(-3t^2u)\,x_b(-3tu^2),
\end{equation}

\begin{equation}\label{pf27}
[x_{-a-b}(u),x_{2a+b}(t)]=
 x_a(-2tu)\,x_{3a+b}(-3t^2u)\,x_{-b}(3tu^2),
\end{equation}

\begin{equation}\label{pf28}
[x_{-3a-b}(u),x_{2a+b}(t)]
 =x_{-a}(-tu)\,x_{a+b}(-t^2u)\,x_{3a+2b}(t^3u)\,x_b(-t^3u^2),
\end{equation}

\begin{equation}\label{pf29}
[x_{2a+b}(u),x_{-3a-b}(t)]
 =x_{-a}(tu)\,x_{a+b}(tu^2)\,x_{3a+2b}(-tu^3)\,x_b(-2t^2u^3),
\end{equation}

\begin{equation}\label{pf30}
[x_{-3a-2b}(u),x_{2a+b}(t)]
 =x_{-a-b}(tu)\,x_a(-t^2u)\,x_{3a+b}(-t^3u)\,x_{-b}(t^3u^2),
\end{equation}

\begin{equation}\label{pf31}
[x_{2a+b}(u),x_{-3a-2b}(t)]
 =x_{-a-b}(-tu)\,x_a(tu^2)\,x_{3a+b}(tu^3)\, x_{-b}(2t^2u^3),
\end{equation}

\begin{equation}\label{pf32}
[x_{-a}(u),x_{3a+b}(t)]
 =x_{2a+b}(tu)\,x_{a+b}(-tu^2)\,x_b(tu^3)\, x_{3a+2b}(2t^2u^3),
\end{equation}

\begin{equation}\label{pf33}
[x_{3a+b}(u),x_{-a}(t)]
 =x_{2a+b}(-tu)\,x_{a+b}(t^2u)\,x_b(-t^3u)\,x_{3a+2b}(t^3u^2),
\end{equation}

\begin{equation}\label{pf34}
[x_{-2a-b}(u),x_{3a+b}(t)]
 =x_a(-tu)\,x_{-a-b}(tu^2)\,x_{-3a-2b}(tu^3)\, x_{-b}(-2t^2u^3),
\end{equation}

\begin{equation}\label{pf35}
[x_{3a+b}(u),x_{-2a-b}(t)]
 =x_a(tu)\,x_{-a-b}(-t^2u)\,x_{-3a-2b}(-t^3u)\,x_{-b}(-t^3u^2),
\end{equation}

\begin{equation}\label{pf36}
[x_{-3a-2b}(u),x_{3a+b}(t)]=x_{-b}(-tu),
\end{equation}

\begin{equation}\label{pf37}
[x_{-a-b}(u),x_{3a+2b}(t)]
 =x_{2a+b}(-tu)\,x_a(-tu^2)\,x_{-b}(tu^3)\, x_{3a+b}(2t^2u^3),
\end{equation}

\begin{equation}\label{pf38}
[x_{3a+2b}(u),x_{-a-b}(t)]
 =x_{2a+b}(tu)\,x_a(t^2u)\,x_{-b}(-t^3u)\,x_{3a+b}(t^3u^2),
\end{equation}

\begin{equation}\label{pf39}
[x_{-2a-b}(u),x_{3a+2b}(t)]
 =x_{a+b}(tu)\,x_{-a}(tu^2)\,x_{-3a-b}(-tu^3)\, x_b(2t^2u^3),
\end{equation}

\begin{equation}\label{pf40}
[x_{3a+2b}(u),x_{-2a-b}(t)]
 =x_{a+b}(-tu)\,x_{-a}(-t^2u)\,x_{-3a-b}(t^3u)\,x_b(t^3u^2),
\end{equation}

\begin{equation}\label{pf41}
[x_{-3a-b}(u),x_{3a+2b}(t)]=x_{a+b}(-tu).
\end{equation}

\bigskip

Sergey G. Kolesnikov

Siberian Federal University

email: skolesnikov@sfu-kras.ru

\end{document}